\documentstyle[amssymb]{bcp98e}
\newcounter{defn}
\newtheorem{dfn}[defn]{Definition}
\newtheorem{theorem}{Theorem}

\newtheorem{corollary}[theorem]{Corollary}
\newtheorem{lemma}[theorem]{Lemma}

\mathchardef\za="710B  
\mathchardef\zb="710C  
\mathchardef\zg="710D  
\mathchardef\zd="710E  
\mathchardef\zve="710F 
\mathchardef\zz="7110  
\mathchardef\zh="7111  
\mathchardef\zvy="7112 
\mathchardef\zi="7113  
\mathchardef\zk="7114  
\mathchardef\zl="7115  
\mathchardef\zm="7116  
\mathchardef\zn="7117  
\mathchardef\zx="7118  
\mathchardef\zp="7119  
\mathchardef\zr="711A  
\mathchardef\zs="711B  
\mathchardef\zt="711C  
\mathchardef\zu="711D  
\mathchardef\zvf="711E 
\mathchardef\zq="711F  
\mathchardef\zc="7120  
\mathchardef\zw="7121  
\mathchardef\ze="7122  
\mathchardef\zy="7123  
\mathchardef\zf="7124  
\mathchardef\zvr="7125 
\mathchardef\zvs="7126 
\mathchardef\zf="7127  
\mathchardef\zG="7000  
\mathchardef\zD="7001  
\mathchardef\zY="7002  
\mathchardef\zL="7003  
\mathchardef\zX="7004  
\mathchardef\zP="7005  
\mathchardef\zS="7006  
\mathchardef\zU="7007  
\mathchardef\zV="7008  
\mathchardef\zW="700A  

\newcommand{\be}{\begin{equation}}
\newcommand{\ee}{\end{equation}}
\newcommand{\ra}{\rightarrow}

\newcommand{\bea}{\begin{eqnarray}}
\newcommand{\eea}{\end{eqnarray}}
\newcommand{\beas}{\begin{eqnarray*}}
\newcommand{\eeas}{\end{eqnarray*}}

\newcommand{\R}{{\Bbb R}}
\newcommand{\1}{{\bold 1}}

\newcommand{\ti}{\times}
\newcommand{\A}{{\cal A}}
\newcommand{\Li}{{\cal L}}
\newcommand{\X}{{\cal X}}
\newcommand{\C}{{\cal C}}
\newcommand{\Ha}{{\cal H}}
\newcommand{\K}{{\cal K}}
\newcommand{\M}{{\cal M}}
\newcommand{\F}{{\cal F}}

\begin{document}
\mathclass{Primary 17B60 17B66; Secondary 53C15.}

\abbrevauthors{J. Grabowski and K. Grabowska}
\abbrevtitle{The Lie algebra of a Lie algebroid}

\title{The Lie algebra of a Lie algebroid}

\author{Janusz\ Grabowski}
\address{Institute of Mathematics, University of Warsaw\\
ul. Banacha 2, 02-097 Warszawa, Poland\\ and\\
Mathematical Institute, Polish Academy of Sciences\\
ul. \'Sniadeckich 8, P.O. Box 137, 00-950 Warszawa, Poland\\
E-mail: jagrab@mimuw.edu.pl}

\author{Katarzyna\ Grabowska}
\address{Division of Mathematical  Methods  in  Physics,  University  of
Warsaw\\ ul. Ho\.za 74, 00-682 Warszawa, Poland\\
E-mail: konieczn@fuw.edu.pl}

\vfootnote{}{Research supported by KBN, grant No. 2 P03A 041 18.}

\maketitlebcp

\abstract{We  present  results  describing  Lie   ideals   and   maximal
finite-codimensional Lie subalgebras of the Lie algebras associated with
Lie algebroids with non-singular anchor maps. We also prove that every
isomorphism  of  such  Lie  algebras  induces  diffeomorphism  of   base
manifolds respecting the generalized foliations defined  by  the  anchor
maps.}

\section{Introduction.}

One of the most striking fact about Lie algebras of vector fields is the
classical result  of Pursell and Shanks \cite{16} which states that the
Lie algebra ${\cal X}_c(M)$ of  all compactly supported smooth  vector
fields  on    a smooth   manifold $M$
determines the smooth structure of  $M$, i.e., the  Lie  algebras
${\cal X}_c(M_1)$ and ${\cal X}_c(M_2)$ are isomorphic if and only
if $M_1$ and $M_2$ are diffeomorphic.
There are similar results in  special geometric situations
(hamiltonian, contact, group invariant vector  fields, etc.),  as
for  example  the results   of  Omori \cite{15} (Chapter X), Abe
\cite{1},  Atkin   and   Grabowski  \cite{3},  or  Hauser  and  M\"uller
\cite{HM}, for which specific tools were developed in each case,
and in the case of Lie algebras of vector fields which are modules over
the corresponding rings of functions  (we  shall  call  them  {\it
modular}), when the answer is more or less complete.
The standard model of a modular Lie  algebra  of vector fields
is  the  Lie   algebra ${\cal X}({\cal  F})$  of  all  vector  fields
tangent   to   a   given (generalized) foliation ${\cal F}$.
Let us recall the work  of Amemiya  \cite{2}, Grabowski \cite{5}, where
an algebraic  approach  made it possible  to consider analytic cases as
well, and finally the brilliant purely algebraic result of Skryabin
\cite{Sk}. This final result  states, in the geometric situation, that in
the case  when modular  Lie algebras  of   vector   fields contain
finite families of vector fields with no common  zeros  (we  shall
say that they are {\it strongly non--singular}), isomorphisms
between  them  are  generated  by  isomorphisms  of  corresponding
algebras  of  functions,  i.e.,  diffeomorphisms   of   underlying
manifolds.
The method of Shanks and Pursell depends  on  the  description  of
maximal ideals in the Lie algebra  ${\cal  X}_c(M)$  in  terms  of
points of $M$ : maximal ideals are of the form $\Li_p$ for  $p\in
M,$ where $\Li_p$ consists of vector fields  which  are  flat  at
$p.$  However,  this  method  fails  in  analytic  cases,    since
analytic  vector  fields  which  are  flat  at  $p$  are  zero  on   the
whole component of $M.$ Therefore in  \cite{2}  and  \cite{5}
maximal  finite codimensional subalgebras are used instead of ideals.

The notion of a {\it Lie algebroid} generalizes simultaneously the notion
of a Lie algebra and  that  of  the  tangent  bundle  over  a  manifold.
Sections   of   the   Lie   algebroid   bundle    form    a    (possibly
infinite-dimensional) Lie algebra with  properties, generically, close
to that of Lie algebras of vector fields.
In this note we  present   a  similar  approach   to  the  Lie  algebras
associated with Lie algebroids (or, in most  of  the  abstract  setting,
their  algebraic  counterpart:  Lie   pseudoalgebras),   especially   to
describe ideals  and maximal finite-codimensional subalgebras together with
some remarks  concerning Shanks-Pursell type theorems.
Of course, since the theory of  Lie  algebroids includes the theory of
finite-dimensional  Lie  algebras  for  which  no general result is
possible, we impose  some  non-singularity  conditions on the anchor
map.

\section{Lie algebroids and Lie pseudoalgebras.}

Let us start with fixing the terminology and notation.

The geometric part of our results will include as well smooth as
real-analytic and holomorphic cases. Therefore we shall deal at the same
time with finite  dimensional  manifolds  $M$ and vector bundles over $M$
of   different   classes   of smoothness:
${\cal C}=C^\infty  ,  C^\omega   ,   {\cal H}$,
where $C^\infty$ denotes the classical smooth case,  $C^\omega$  -
the real analytic case, and ${\cal H}$  denotes the holomorphic  case
for Stein manifolds. For details we refer to \cite{3}.
For instance, ${\cal C}(M)$ is  the  algebra  of  class ${\cal  C}$
functions on the manifold $M$ of class ${\cal C}$.  Note  that  the
algebras $C^\infty (M)$  and  $C^\omega  (M)$  are  real  and  the
algebra ${\cal  H}(M)$  of  holomorphic  functions  on  the  Stein
manifold $M$ is complex. It is well known  that the  corresponding
Lie algebra ${\cal X}(M)$ of all class ${\cal  C}$  vector  fields
can be regarded as the Lie algebra of derivations of ${\cal C}(M)$
(in the analytic cases we refer to \cite{6}).

\begin{dfn} Let $M$ be a manifold of class $\C$. A {\rm Lie algebroid  on
$M$} is a vector bundle $\zt:L\ra M$, together with a  bracket
$[\cdot,\cdot]:\zG L\ti \zG L\ra\zG L$ on the module $\zG L$ of class
$\C$ global sections of $L$, and a vector bundle morphism $a:L\ra TM$,
over  the  identity in $M$, from $L$ to the tangent bundle $TM$, called
the {\rm  anchor}  of the Lie algebroid, such that
\noindent
\item{(i)} the bracket on $\zG L$ is $\R$-bilinear (or $\Bbb C$-bilinear
in the case $\C=\Ha$), alternating, and satisfies the Jacobi identity;
\item{(ii)} $[X,fY]=f[X,Y]+a(X)(f)Y$ for all  $X,Y\in  \zG  L$  and  all
$f\in\C(M)$;
\item{(iii)} $a([X,Y])=[a(X),a(Y)]$ for all $X,Y\in\zG L$.
\end{dfn}

We get an algebraic counterpart of the notion of Lie algebroid replacing
$\C(M)$ by an arbitrary algebra $\A$, and the  sections  of  the  vector
bundle $\zt:L\ra M$ by a module $\Li$ over the algebra $\A$.

\begin{dfn} Let $R$ be a commutative and unitary ring, and let $\A$ be a
commutative and unitary $R$-algebra. A {\rm Lie pseudoalgebra  over  $R$
and $\A$} is an $\A$-module $\Li$ together with a bracket $[\cdot,\cdot]
:\Li\ti\Li\ra\Li$ on the module $\Li$, and an $\A$-module  morphism
$a:\Li\ra{\rm Der}(\A)$ from $\Li$ to the $\A$-module ${\rm Der}(\A)$ of
derivations of $\A$, called the {\rm anchor} of $\Li$, such that
\noindent
\item{(i)} the bracket on $\Li$ is $R$-bilinear , alternating, and
satisfies the Jacobi identity;
\item{(ii)} For all  $X,Y\in  \Li$  and  all $f\in\A$ we have
\be\label{0}
[X,fY]=f[X,Y]+a(X)(f)Y;
\ee
\item{(iii)} $a([X,Y])=[a(X),a(Y)]$ for all $X,Y\in\Li$.
\end{dfn}
Lie algebroids on a singleton base space are Lie algebras. Another
extreme example is the tangent bundle $TM$ with the canonical bracket
on the space $\X=\zG TM$ of vector fields.

Lie pseudoalgebras appeared first in the paper of Herz  \cite{He}
but one can find similar concepts under more than a dozen  of  names  in
the   literature   (e.g.   $(R,A)$-Lie   algebras,   Lie-Cartan   pairs,
Lie-Rinehart algebras, differential algebras, etc.). Lie algebroids were
introduced  by Pradines \cite{Pr}. For both notions we refer to a
survey article  by Mackenzie \cite{Ma}.

A subset $V$ in an $\A$-module $E$ which is  an  $\A$-submodule  in  $E$
will be called {\it modular}.
When the Lie pseudoalgebra (so its anchor) is fixed, we shall write  shortly
$\hat X$ instead of  $a(X)$.  Thus  $\hat\Li=\{\hat  X:X\in\Li\}$  is  a
modular  Lie subalgebra  of  ${\rm  Der}(\A)$.  We shall  call  a  Lie
pseudoalgebra $\Li$  {\it  strongly  non-singular}  if  $\hat\Li$  is  a
strongly  non-singular  Lie subalgebra  in  ${\rm   Der}(\A)$   in   the
terminology  of  \cite{9},  i.e.,  $\hat\Li(\A)=\A$  with  the  obvious
notation $\hat\Li(\A)=span\{\hat X(f):X\in\Li,\ f\in\A\}$.

We  shall  be  interested  in  the  Lie  algebra  structure  of   a   Lie
pseudoalgebra $\Li$, i.e., we shall consider $\Li$ as a Lie algebra  over
$R$.

For $I\subset\A$ we denote $\Li_I=\{ X\in\Li:\hat X(\A)\subset I\}$.  It
is obvious that if $I$ is an ideal of $\A$, then $I\Li\subset\Li_I$  and
$I\Li, \Li_I$ are modular Lie subalgebras of  $\Li$.  If,  additionally,
$I$ is $\Li$-invariant, i.e., $\hat\Li(I)\subset  I$,  then  $I\Li$  and
$\Li_I$ are Lie ideals of $\Li$. The kernel $\Li_0$ of the anchor map is
also a modular Lie ideal of $\Li$. The quotient  algebra  $\Li/\Li_0$  is
canonically isomorphic with the modular Lie subalgebra $\hat\Li=a(\Li)$ in
${\rm Der}(\A)$. Note that in the Lie algebroid case $\hat\Li$ generates
a generalized distribution which is involutive due to $(iii)$ of the
definition. This distribution is finitely generated by vector fields
from $\hat\Li$. Indeed, it is well known that the module $\Li$ of
sections of a vector bundle $L$ is finitely generated (and projective), say
by $X_1,\dots,X_k$, so that $a(X_1),\dots,a(X_k)$ generate $\hat\Li$. It
is also well known that in this case the generalized distribution is
completely integrable and defines a generalized foliation in the sense of
Sussman \cite{Su} and Stefan \cite{St}.

\section{Lie ideals.}

Throughout this section we fix a Lie pseudoalgebra over $R$ and $\A$ and
assume that 2 is invertible in the ring $R$.

The following Lemma is essentially due to Skryabin \cite{Sk}  (cf.
also \cite{8}).

\begin{lemma} If $\Li$ is strongly non-singular, then
\be
\hat X(\A)\Li\subset[\Li,[\Li,X]]
\ee
for every $X\in\Li$.
\end{lemma}
\Proof Take $Y\in\Li$. For $f,g\in\A$ set $B(f,g)=[fY,[gY,X]]$. It is  a
matter  of  simple  calculations, using  (\ref{0}), to show the identity
\be
B(fg^2,1)-2B(fg,g)+B(f,g^2)=-4f\hat Y(g)\hat X(g)Y,
\ee
so
\be\label{1}
f\hat Y(g)\hat X(g)Y\in[\Li,[\Li,X]].
\ee
After the linearization  with
respect to $Y$ (i.e., we set $Y:=Y+Z$ in (\ref{1}) and use (\ref{1})  to
remove the terms not mixing $Y$ and $Z$), we conclude that
\be\label{2}
f\hat Y(g)\hat X(g)Z+f\hat Z(g)\hat X(g)Y\in[\Li,[\Li,X]],
\ee
for all $f,g\in\A$  and  all  $Y,Z\in\Li$.  Putting  $f:=\hat  Y(g)$  in
(\ref{2}) and using (\ref{1}), we get
\be
(\hat Y(g))^2\hat X(g)Z\in[\Li,[\Li,X]]
\ee
and further (after $g:=f+g$, $Z:=\hat Y(g)Z$), in a similar way,
\be\label{3}
(\hat Y(g))^3\hat X(f)Z\in[\Li,[\Li,X]],
\ee
for all $f,g\in\A$ and  $Y,Z\in\Li$.  Now, (\ref{3})  implies  that  the
radical ${\rm rad}(J)$ of  the  largest  ideal  $J$   of   $\A$
such   that   $J\hat X(\A)\Li\subset  [\Li,[\Li,X]]$  includes
$\hat\Li(\A)$   and,   since $\hat\Li(\A)=\A$, we have ${\rm rad}(J)=\A$,
so $J=\A$ ($\A$ is unitary) and the lemma follows.
\endproof
\begin{corollary} For a strongly non-singular Lie  pseudoalgebra   $\Li$,
we have $[\Li,\Li]=\Li$.
\end{corollary}
\Proof According to  Lemma  1,   $\hat\Li(\A)\Li\subset[\Li,[\Li,\Li]]$,
so $\Li\subset [\Li,[\Li,\Li]]\subset[\Li,\Li]$.
\endproof
The following theorem describing Lie ideals is a version of \cite{8} for
Lie pseudoalgebras.
\begin{theorem}  Assume  that  $\Li$  is  a  strongly  non-singular   Lie
pseudoalgebra.Then,
\noindent
\item{(a)} given a Lie ideal $\K$ of $\Li$, there  is  a  $\Li$-invariant
(associative) ideal $I$ of $\A$  such  that  $\hat\K(\A)\subset  I$  and
$I\Li\subset\K\subset\Li_I$;
\item{(b)} given a Lie ideal $\K$  of  $\Li$,  the  ideal  $[\Li,\K]$  is
modular. Moreover, $\K$ is modular if and only if $[\Li,\K]=\K$.
\end{theorem}
\Proof Let $I$ be the largest ideal of $\A$ such  that  $I\Li\subset\K$,
i.e., $I=\{ f\in\A:f\Li\subset\K\}$. According to Lemma 1,
\be\label{4}
\hat\K(\A)\Li\subset[\Li,[\Li,\K]]\subset[\Li,\K]\subset\K,
\ee
so that $\hat\K(\A)\subset I$. Moreover, the property (\ref{0}) implies
\be
\hat\Li(I)\Li\subset([\Li,I\Li]+I\Li)\subset([\Li,\K]+\K)\subset\K,
\ee
thus $\hat\Li(I)\subset I$.

Similarly,
\be
\A[\Li,\K]\subset([\A\Li,\K]+\hat\K(\A)\Li)\subset[\Li,\K],
\ee
which shows that $[\Li,\K]$ is modular. If $\K$ itself is modular, then
\be
\K=\hat\Li(\A)\K\subset([\Li,\A\K]+\A[\Li,\K])\subset[\Li,\K],
\ee
since $[\Li,\K]$ is modular, and $\K=[\Li,\K]$ follows.
\endproof

By a {\it maximal ideal} (with a given property) we mean an ideal  which
is maximal in the  family  of  all  ideals  (with  the  given  property)
which are different from the whole  algebra  (but  including  the   zero
ideal). In this sense, the zero ideal is the only  maximal  ideal  in  a
simple algebra. The same terminology holds for subalgebras.

\begin{corollary} Under the assumptions  of   Theorem 3   there   is
a   bijection   between   maximal   $\Li$-invariant   ideals   $I$    of
the associative algebra $\A$ and maximal Lie ideals of $\Li$ given by
$I\mapsto\Li_I$. Moreover, every Lie ideal of $\Li$ is  contained  in  a
maximal one.
\end{corollary}
\Proof First, we show that the mapping $I\mapsto\Li_I$ is  injective  on
the  set  of  maximal  $\Li$-invariant  ideals  of  $\A$.  We shall  show
even more: $\Li_I\subset\Li_J$ implies $I=J$. Indeed,   if
$\Li_I\subset\Li_J$ for maximal $\Li$-invariant ideals  $I,J$  of  $\A$,
then                   $(I+J)\Li=(I\Li+J\Li)\subset(\Li_I+\Li_J)=\Li_J$,
since $I\Li\subset\Li_I$   and   $J\Li\subset\Li_J$.   But,
if $I\ne J$, then $(I+J)$   is a $\Li$-invariant  ideal  of  $\A$  larger
than  $I$,  so  $(I+J)=A$ and $\A\Li=\Li\subset\Li_J$;  a  contradiction,
since  $\Li$  is strongly non-singular, so $\Li_J\ne \Li$.

Since the algebra $\A$ is unital, the union  of  any  increasing
chain of $\Li$-invariant ideals is again  a  $\Li$-invariant  ideal  and
different from $\A$, so that every $\Li$-invariant ideal is contained in
a maximal $\Li$-invariant ideal of  $\A$.  In  view  of  Theorem 3  any
proper Lie ideal of $\Li$ is contained in $\Li_I$ for a $\Li$-invariant,
hence maximal $\Li$-invariant, ideal $I$ of $\A$. It  remains  to  show
that $\Li_I$ is maximal. Indeed, if $\K$ is a larger proper  Lie  ideal,
then $\K\subset\Li_J$ for a maximal $\Li$-invariant ideal $J$  of  $\A$.
But we already know that this implies $I=J$, so $\Li_I=\Li_J=\K$.
\endproof

\begin{corollary} If $\A$ is a  simple  $\Li$  module,  i.e.,  the  only
$\Li$-invariant ideals are $\A$ and $\{ 0\}$, then $\Li$ is a simple Lie
algebra.
\end{corollary}

In the case of Lie algebroids, the  maximal  $\Li$-invariant  ideals  of
$\C(M)$ are known to consist of functions which are flat at a point when
restricted  to  a  leaf  of  the  generalized  foliation  determined  by
$a(\Li)$; in analytic cases this means that they are zero on the closure
of this leaf. The corresponding Lie ideals consist of  those  $X\in\Li$
which are mapped by the anchor map to vector fields  which  are  flat  at
corresponding points of the leaves (in the analytic cases: vanish on the
corresponding leaves). In particular, the Lie algebra  of  real-analytic
vector fields on a compact real-analytic  manifold  is  simple,  as  was
first proved in \cite{5}.

\section{Finite-codimensional subalgebras and isomorphisms.}

Throughout this section we assume that $R$ is a field of  characteristic
0.

By $\M(\A)$ we denote the set of all maximal finite-codimensional ideals
of $\A$ and by $\M(\Li)$ -- the set of all maximal  finite-codimensional
Lie subalgebras of $\Li$. It is well known  (cf.  \cite{5},  Proposition
3.5) that in the case $\A=\C(M)$ we may identify $\M(\A)$ with $M$ by
\be
M\ni p\mapsto J(p)=\{ f\in\C(M):f(p)=0\}\in\M(\A).
\ee
For $I\subset\A$, let us set $V(I)=\{ J\in\M(\A): I\subset J\}$
and
\be
\bar I=\bigcap_{J\in V(I)}J.
\ee
In \cite{2,5,8} points of the manifold $M$ were represented  by  maximal
finite-codimensional Lie subalgebras of the corresponding  Lie  algebras
of vector fields.  Then,  this  description  was  used for associating a
diffeomorphism with a given isomorphism of such Lie algebras. A  similar
result can be proved for strongly non-singular Lie pseudoalgebras.

\begin{theorem} Let $\Li$ be a strongly non-singular  Lie  pseudoalgebra
over  $R$  and   $\A$.   Then,   for    any   finite-codimensional   Lie
subalgebra  $\K$  of $\Li$ there is a finite-codimensional  ideal  $I$
of  $\A$ such that $I\Li\subset\K\subset\Li_{\bar I}$.
\end{theorem}

\Proof  Since  $\hat\Li(\A)=\A$, there are $X_1,\dots,X_m\in\Li$ and
$f_1,\dots,f_m\in\A$  such  that  $\sum_{i=1}^m\hat  X_i(f_i)=\1$. It is
easy   to   see   that   $\K_1=\{   X\in\K:[\Li,X]\subset\K\}$   is    a
finite-codimensional Lie subalgebra  of  $\Li$  as  the  kernel  of  the
adjoint representation of $\K$ in the finite-dimensional space $\Li/\K$.
Similarly,     $\K_2=\{     X\in\K_1:[\Li,X]\subset\K_1\}$     is      a
finite-codimensional    Lie    subalgebra    of     $\Li$.     Moreover,
$[\Li,[\Li,\K_2]]\subset\K$. Hence,
\be
W=\{ f\in\A: fX_i\in\K_2,\ i=1,\dots,m\}
\ee
is a finite-codimensional subspace of $\A$ and, in view of Lemma 1,
\be
W\hat X_i(f_i)\Li\subset[\Li,[\Li,K_2]]\subset\K.
\ee
Hence
\be
W\left(\sum_{i=1}^m\hat X_i(f_i)\right)\Li\subset\K,
\ee
so that $I\Li\subset\K$, where $I$ is the ideal  of  $\A$  generated  by
$W$, thus finite-codimensional. We can assume that $I$ is  maximal  with
this property. Then, according to (\ref{0}),
\be
\hat\K(I)\Li\subset([\K,I\Li]+I[\K,\Li])\subset\K,
\ee
so    that    $\hat\K(I)\subset    I$.    This    in    turn     implies
$\hat\K(\A)\subset\bar I$, as show Lemmata 4.1 and 4.2 in \cite{5}.
\endproof

\begin{corollary} If $\Li_J$ is of finite codimension in $\Li$ for each
$J\in\M(\A)$, then the map $J\mapsto\Li_J$  constitutes  a  bijection  of
$\M(\A)$ with $\M(\Li)$. In particular, for $\Li$ being a Lie  algebroid
of class $\C$ on $M$, we have the bijection
\be
M\ni p\mapsto\Li_p=\{ X\in\Li: \hat X(p)=0\}\in\M(\Li).
\ee
\end{corollary}
\Proof The proof is parallel to that of Corollary 4.
\endproof
Now, assume that $\Li$ is a Lie algebroid on $M$. It is clear that
the kernel $\Li_0$ of the anchor map equals $\bigcap_{p\in M}\Li_p$, so
that it can be defined in terms of the Lie algebra structure of $\Li$.
Namely, the kernel of the anchor map is the intersection
\be\label{an}
\Li_0=\bigcap_{\K\in\M(\Li)}\K.
\ee
If $\zV:\Li\ra\Li$ is an automorphism of the Lie algebra $\Li$, then $\zV$
maps maximal finite-codimensional subalgebras into maximal
finite-codimensional subalgebras, so that it preserves their intersection
and, in view (\ref{an}),  $\zV(\Li_0)=\Li_0$. This, in turn, implies
that $\zV$ induces an automorphism $\hat\zV:\hat\Li\ra\hat\Li$ of the Lie
algebra $\hat\Li=\Li/\Li_0$ of vector fields on $M$. Since
$\hat\Li(\C(M))=\C(M)$, known results on isomorphism of Lie algebras of
vector fields (e.g. \cite{5}, Theorem 5.5, or  \cite{Sk}, Theorem 3.2)
imply that $\hat\zV$ is generated by an automorphism  of $\C(M)$, i.a., a
class $\C$ diffeomorphism of $M$.
Of course, a similar reasoning remains valid for two Lie algebroids with
the required properties of their anchor maps and we get the following.

\begin{theorem} Let $\Li^i$ be the Lie algebra of sections of a strongly
non-singular Lie algebroid of class $\C$ on a manifold $M^i$, with an
anchor map $a^i:\Li^i\ra\X(M^i)$ and the generalized foliation $\F^i$
generated by $\hat\Li^i=a^i(\Li^i)$, $i=1,2$.
If $\zV:\Li^1\ra\Li^2$ is an isomorphism of the Lie algebras $\Li^1$ and
$\Li^2$, then $\zV$ maps the kernel $\Li_0^1$ of the anchor map $a^1$ onto
the kernel $\Li_0^2$ of the anchor map $a^2$ and therefore it
induces an isomorphism $\hat\zV:\hat\Li^1\ra\hat\Li^2$ of the
corresponding Lie algebras of vector fields. Moreover, $\hat\zV$ is
implemented by a class $\C$ diffeomorphism $\phi:M^1\ra M^2$ (i.e.,
$\hat\zV=\phi_*$) which preserves the generalized foliations:
$\phi(\F^1)=\F^2$.
\end{theorem}

{\nine

}
\end{document}